# A kernel type nonparametric density estimator for decompounding

BERT VAN ES[*], SHOTA GUGUSHVILI[**] and PETER SPREIJ[†]

*Korteweg-de Vries Institute for Mathematics, Universiteit van Amsterdam, Plantage Muidergracht 24, 1018 TV Amsterdam, The Netherlands. E-mail: [*]vanes@science.uva.nl; [**]sgugushv@science.uva.nl; [†]spreij@science.uva.nl*

Given a sample from a discretely observed compound Poisson process, we consider estimation of the density of the jump sizes. We propose a kernel type nonparametric density estimator and study its asymptotic properties. An order bound for the bias and an asymptotic expansion of the variance of the estimator are given. Pointwise weak consistency and asymptotic normality are established. The results show that, asymptotically, the estimator behaves very much like an ordinary kernel estimator.

*Keywords:* asymptotic normality; consistency; decompounding; kernel estimation

## 1. Introduction

Let $N(\lambda)$ be a Poisson random variable with parameter $\lambda$ and let $Y_1, Y_2, \ldots$ be a sequence of independent and identically distributed random variables that are independent of $N(\lambda)$, have a common distribution function $F$ and have density $f$. Consider a Poisson sum of $Y$s:

$$X = \sum_{j=1}^{N(\lambda)} Y_j.$$

Assume $\lambda$ is known. The statistical problem we consider is nonparametric estimation of the density $f$ based on observations on $X$. Because adding a Poisson number of $Y$s is referred to as compounding, we refer to the problem of recovering the density $f$ of $Y$s from the observations on $X$ as decompounding. The problem of estimating the density $f$ is equivalent to the problem of estimating the jump size density $f$ of a compound Poisson process $X' = (X'_t)_{t \geq 0}$ with intensity $\lambda$ when the process is observed at equidistant time points (rescaling if necessary, the observation step size can be taken to be equal to 1). Compound Poisson processes have important applications in queueing and risk theory (see, e.g., Embrechts *et al.* [7] and Prabhu [11]), for example, the random variables $Y_1, Y_2, Y_3, \ldots$ can be interpreted as claims of random size that arrive at an insurance company or as the number of customers who arrive at a service point at random times with exponentially distributed interarrival time.

The problem of nonparametric estimation of the distribution function $F$ in the case of both continuous and discrete laws was treated by Buchmann and Grübel [1]. Their estimation method is based on a suitable inversion of the compounding operation (i.e., transition from the distribution of $Y$ to the distribution of $X$) and use of an empirical estimator for the distribution of $X$, thus resulting in a plug-in type estimator for the distribution of $Y$. A further ramification of this approach in the case of a discrete law was given by Buchmann and Grübel [2]. To the best of our knowledge, the present paper is the first attempt to (nonparametrically) estimate the density $f$. A very natural use of nonparametric density estimators is informal investigation of the properties of a given set of data. The estimators can give valuable indications about the shape of the density function, for example, such features as skewness and multimodality. The knowledge of these features might come in handy in applications, for example, in insurance, where $f$ is a claim size density.

One possible way to construct an estimator for the density $f$ (suggested in Hansen and Pitts [9]) is via smoothing the plug-in type estimator $F_n$ of the distribution function $F$, that was defined by Buchmann and Grübel [1], with a kernel, but at present no theoretical results for this estimator seem to be available. We opt for an alternative approach based on inversion of the characteristic function $\phi_f$, an approach that is in spirit similar to the use of kernel estimators in deconvolution problems (the latter were first introduced by Liu and Taylor [10] and Stefanski and Caroll [15]; for a more recent overview, see Wand and Jones [19]). Before we proceed any further, we need to specify the observation scheme. Zero observations provide no information on the $Y$s and, hence, an estimator of $f$ should be based on nonzero observations. In a sample of fixed size there are a random number of nonzero observations. We want to avoid this extra technical complication, so we assume that we have observations $X_1, \ldots, X_{T_n}$ on $X$, where $T_n$ is the first moment we get precisely $n$ nonzero observations ($T_n$ of course is random). We denote the nonzero observations by $Z_1, Z_2, \ldots, Z_n$.

We turn to the construction of the estimator of the density $f$. First note that the characteristic function of $X$ is given by

$$\phi_X(t) = \mathrm{E}[\mathrm{e}^{\mathrm{i}tX}] = \mathrm{e}^{-\lambda + \lambda \phi_f(t)},$$

where $\phi_f$ denotes the characteristic function of a random variable with density $f$. Rewrite the characteristic function of $X$ as

$$\phi_X(t) = \mathrm{e}^{-\lambda} + (1 - \mathrm{e}^{-\lambda}) \frac{1}{\mathrm{e}^\lambda - 1}(\mathrm{e}^{\lambda \phi_f(t)} - 1).$$

Denote the density of $X$ given $N > 0$ by $g$. It follows that the characteristic function of $X$ given $N > 0$ is equal to

$$\phi_g(t) = \frac{1}{\mathrm{e}^\lambda - 1}(\mathrm{e}^{\lambda \phi_f(t)} - 1).$$

Because $\phi_f$ vanishes at plus and minus infinity, so does $\phi_g$. By inverting the above relationship, we get

$$\phi_f(t) = \frac{1}{\lambda}\mathrm{Log}((\mathrm{e}^\lambda - 1)\phi_g(t) + 1).$$

Here Log denotes the *distinguished logarithm* (in general, we cannot use a principal branch of the logarithm) and we refer to Chung ([4], Theorem 7.6.2) and Finkelestein *et al.* [8] for details of its construction. Notice that whenever $\lambda < \log 2$, the distinguished logarithm reduces to the principal branch of an ordinary logarithm. By Fourier inversion, if $\phi_f$ is integrable, we have

$$f(x) = \frac{1}{2\pi\lambda}\int_{-\infty}^{\infty}\mathrm{e}^{-\mathrm{i}tx}\mathrm{Log}((\mathrm{e}^\lambda - 1)\phi_g(t) + 1)\,\mathrm{d}t. \tag{1.1}$$

This relation suggests that if we construct an estimator of $g$ (and hence of $\phi_g$), we will automatically get an estimator for $f$ by a plug-in device. Let $w$ denote a kernel function with characteristic function $\phi_w$ and let $h$ denote a positive number—the bandwidth. The density $g$ will be estimated by the kernel density estimator

$$g_{nh}(x) = \frac{1}{n}\sum_{j=1}^{n}\frac{1}{h}w\biggl(\frac{x - Z_j}{h}\biggr).$$

Properties of kernel estimators can be found in recent books such as Devroye and Györfi [6], Prakasa Rao [12], Tsybakov [16] and Wand and Jones [19]. The characteristic function $\phi_{g_{nh}}$ serves as an estimator of $\phi_g$ and is equal to $\phi_{\mathrm{emp}}(t)\phi_w(ht)$, where $\phi_{\mathrm{emp}}$ denotes the empirical characteristic function

$$\phi_{\mathrm{emp}}(t) = \frac{1}{n}\sum_{j=1}^{n}\mathrm{e}^{\mathrm{i}tZ_j}.$$

In view of (1.1) it is tempting to introduce an estimator

$$\frac{1}{2\pi\lambda}\int_{-\infty}^{\infty}\mathrm{e}^{-\mathrm{i}tx}\mathrm{Log}((\mathrm{e}^\lambda - 1)\phi_{\mathrm{emp}}(t)\phi_w(ht) + 1)\,\mathrm{d}t, \tag{1.2}$$

but there are two problems. First, the measure of those $\omega$s from the underlying sample space $\Omega$ for which the path $(\mathrm{e}^\lambda - 1)\phi_{g_{nh}}(t) + 1$ can become zero is positive (although as $n \to \infty$, this probability tends to zero) and the distinguished logarithm cannot be defined for such $\omega$s. Second, there is no guarantee that the integral in (1.2) is finite. Therefore, we will make the adjustments

$$\hat{f}_{nh}(x) = (M_n \wedge f_{nh}(x)) \vee (-M_n), \tag{1.3}$$

where for those $\omega$s for which the paths $(e^\lambda - 1)\phi_{\text{emp}}(t)\phi_w(ht) + 1$ do not vanish, $f_{nh}$ is given by

$$f_{nh}(x) = \frac{1}{2\pi\lambda} \int_{-1/h}^{1/h} e^{-itx} \text{Log}((e^\lambda - 1)\phi_{\text{emp}}(t)\phi_w(ht) + 1) \, dt$$

and is zero otherwise. Here $M = (M_n)_{n \geq 1}$ is a sequence of positive real numbers that converge to infinity at a suitable rate. We also assume that $\phi_w$ is supported on $[-1, 1]$. Of course, for the truncation in (1.3) to make sense, $f_{nh}(x)$ must be real-valued, but this is easy to check through the change of the integration variable from $t$ into $-t$.

The rest of the paper is organized as follows. Section 2 contains the main results of the paper. In it we derive an order bound for the bias and an asymptotic expansion of the variance of $\hat{f}_{nh}$ at a fixed point $x$, and we show that the estimator is weakly consistent and asymptotically normal. Section 3 provides some simulation results. All the proofs are collected in Section 4.

## 2. Asymptotic properties of the estimator

As is usual in nonparametric estimation, the nonparametric setting forces us to make some smoothness assumptions on the density $f$. Let $\beta, L_1$ and $L_2$ denote some positive numbers and let $l = \lfloor \beta \rfloor$ denote the integer part of $\beta$. If $l = 0$, then by definition set $f^{(l)} = f$. Recall the definition of Hölder and Nikol'ski classes of the functions (cf. Tsybakov [16], pages 5, 19).

**Definition 2.1.** *A function $f$ is said to belong to the Hölder class $\mathcal{H}(\beta, L_1)$ if its derivatives up to order $l$ exist and verify the condition*

$$|f^{(l)}(x+t) - f^{(l)}(x)| \leq L_1 |t|^{\beta - l}.$$

**Definition 2.2.** *A function $f$ is said to belong to the Nikol'ski class $\mathcal{N}(\beta, L_2)$ if its derivatives up to order $l$ exist and verify the condition*

$$\left[ \int_{-\infty}^{\infty} (f^{(l)}(x+t) - f^{(l)}(x))^2 \, dx \right]^{1/2} \leq L_2 |t|^{\beta - l}.$$

We formulate the condition on the density $f$.

**Condition F.** *The density $f$ belongs to $\mathcal{H}(\beta, L_1) \cap \mathcal{N}(\beta, L_2)$. Moreover, $t^\beta \phi_f$ is integrable and the derivatives $f', \ldots, f^{(l)}$ are integrable.*

The following lemma holds true. It is proved in Section 4.

**Lemma 2.1.** *Assume that Condition F holds. Then the density $g$ belongs to $\mathcal{H}(\beta, L_1) \cap \mathcal{N}(\beta, \lambda e^\lambda (e^\lambda - 1)^{-1} L_2)$. Moreover, $t^\beta \phi_g(t)$ is integrable.*

We will use this fact in the proofs of Propositions 2.1 and 2.2 and Theorems 2.1 and 2.2. The requirement that $g \in \mathcal{N}(\beta, \lambda e^{\lambda}(e^{\lambda}-1)^{-1}L_2)$ is motivated by the fact that in the proofs we will make use of the expansion of the mean integrated squared error of a kernel density estimator $g_{nh}$ (cf. Tsybakov [16], page 21), while $g \in \mathcal{H}(\beta, L_1)$ is a standard condition in ordinary kernel density estimation (see Tsybakov [16], Proposition 1.2). The integrability of $f', \ldots, f^{(l)}$ is used in the proof of Lemma 2.1.

**Definition 2.3.** *A function $w$ is called a kernel of order $l$ if the functions $u^j w(u)$, $j = 0, \ldots, l$, are integrable and verify the condition*

$$\int_{-\infty}^{\infty} w(u)\,du = 1, \qquad \int_{-\infty}^{\infty} u^j w(u)\,du = 0 \qquad \text{for } j = 1, \ldots, l-1.$$

Because it is generally recognized that the choice of a kernel is less important for the performance of an estimator (see Wand and Jones [19], page 31), we feel free to impose the following condition on the kernel.

**Condition W.** *The kernel function $w$ satisfies the following conditions:*

1. *$w$ is a bounded symmetric kernel of order $l$.*
2. *The characteristic function $\phi_w$ has a support on $[-1, 1]$.*
3. *$\int_{-\infty}^{\infty} |u|^{\beta} |w(u)|\,du < \infty$.*
4. *$\lim_{|u| \to \infty} |u w(u)| = 0$.*
5. *$\phi_w$ is continuously differentiable.*

To get a consistent estimator, we need to control the bandwidth, so we impose the following restriction.

**Condition H.** *The bandwidth $h$ depends on $n$ and is of the form $h = Cn^{-\gamma}$ for $0 < \gamma < 1$, where $C$ is some constant.*

We also formulate the condition on the truncating sequence $M = (M_n)_{n \geq 1}$ (see Section 1).

**Condition M.** *The truncating sequence $M = (M_n)_{n \geq 1}$ is given by $M_n = n^{\alpha}$, where $\alpha$ is some strictly positive number.*

As the performance criterion, we select the mean squared error

$$\mathrm{MSE}[\hat{f}_{nh}(x)] = \mathrm{E}[(\hat{f}_{nh}(x) - f(x))^2].$$

By standard properties of mean and variance

$$\mathrm{MSE}[\hat{f}_{nh}(x)] = (\mathrm{E}[\hat{f}_{nh}(x)] - f(x))^2 + \mathrm{Var}[\hat{f}_{nh}(x)],$$

the sum of the squared bias and variance at $x$.

First we study the behaviour of the bias of the estimator $\hat{f}_{nh}(x)$.

**Proposition 2.1.** *Suppose Conditions F, W, H and M are satisfied. Then the bias of the estimator $\hat{f}_{nh}(x)$ admits an order bound*

$$\mathrm{E}[\hat{f}_{nh}(x)] - f(x) = O\left(h^\beta + \frac{1}{nh}\right).$$

In ordinary kernel estimation, under the assumption $g \in \mathcal{H}(\beta, L_1)$, the bias is of order $h^\beta$ (see Tsybakov [16], Proposition 1.2). We have an additional term of order $(nh)^{-1}$ that comes from the difficulty of the decompounding problem. Under standard conditions $h \to 0$ and $nh \to \infty$, the bias will asymptotically vanish.

*Remark 2.1.* If $\beta = 2$, then as in our technical report [17], it is possible to derive an exact asymptotic expansion for the bias. The leading term in bias expansion will be

$$-h^2 \frac{\sigma^2(e^\lambda - 1)}{4\pi\lambda} \int_{-\infty}^\infty e^{-itx} \frac{t^2 \phi_g(t)}{(e^\lambda - 1)\phi_g(t) + 1} \, dt.$$

Now let us study the variance of the estimator $\hat{f}_{nh}(x)$.

**Proposition 2.2.** *Suppose that apart from Conditions F, W, H and M, an additional condition $nh^{1+4\beta} \to 0$ holds true. Then the variance of the estimator $\hat{f}_{nh}(x)$ admits the decomposition*

$$\mathrm{Var}[\hat{f}_{nh}(x)] = \frac{1}{nh} \frac{(e^\lambda - 1)^2}{\lambda^2} g(x) \int_{-\infty}^\infty (w(u))^2 \, du + o\left(\frac{1}{nh}\right). \tag{2.1}$$

We see that the variance of our estimator is of the same order as the variance of an ordinary kernel estimator (cf. Tsybakov [16], Proposition 1.4). Under the standard assumption $nh \to \infty$, it will vanish. From a practical point of view, the restriction $nh^{1+4\beta} \to 0$ is not restrictive, especially in view of Proposition 2.3 given below.

By combining Propositions 2.1 and 2.2, we get the following corollary.

**Corollary 2.1.** *Suppose Conditions F, W, H and M hold. The estimator $\hat{f}_{nh}(x)$ is pointwise weakly consistent under the additional assumption $nh^{1+4\beta} \to 0$.*

Recall that the bandwidth $h_{\mathrm{opt}}$ that asymptotically minimizes the mean squared error of a kernel estimator is called optimal. From Propositions 2.1 and 2.2 it is now possible to determine the order of the optimal bandwidth for the estimator $\hat{f}_{nh}$.

**Proposition 2.3.** *The optimal bandwidth $h_{\mathrm{opt}}$ is of order $n^{-1/(2\beta+1)}$. Furthermore, the mean squared error of the estimator $\hat{f}_{nh}$ computed for the optimal bandwidth is of order $n^{-2\beta/(2\beta+1)}$.*

Note that the optimal bandwidth is of order $n^{-1/(2\beta+1)}$, just as in the case of ordinary kernel estimation.

**Remark 2.2.** When $\beta = 2$, then as in van Es *et al.* ([17], Proposition 3.3), it is possible to derive an exact expression for $h_{\mathrm{opt}}$:

$$h_{\mathrm{opt}} = \left( \frac{4\pi^2 g(x) \int_{-\infty}^{\infty} (w(u))^2 \, \mathrm{d}u}{\sigma^4 (\int_{-\infty}^{\infty} \mathrm{e}^{-\mathrm{i}tx} t^2 \phi_g(t)/((\mathrm{e}^{\lambda}-1)\phi_g(t)+1) \, \mathrm{d}t)^2} \right)^{1/5} n^{-1/5}.$$

The extension of our results to the data-dependent bandwidth case is outside the scope of the present paper.

It is interesting to verify whether our estimator is minimax. We refer to van Es *et al.* ([17], Theorem 3.1), where for $\beta = 2$ we proved that the minimax convergence rate for a quadratic loss function is at least $n^{-2/5}$ and, that our estimator attains it for a fixed density $f$. This result can be easily generalized to an arbitrary $\beta > 0$. Whether the estimator itself is minimax is an open question. In any case, the results of the present section show that its behaviour is rather reasonable.

Concluding this section, we will derive two asymptotic normality results for $\hat{f}_{nh}$.

**Theorem 2.1.** *Assume that the Conditions F, W, H and M hold, and that the bandwidth $h$ satisfies an additional condition $nh^{2\beta+1} \to 0$ and $g(x) \neq 0$. Then*

$$\left( \frac{\hat{f}_{nh}(x) - f(x)}{\sqrt{\mathrm{Var}[\hat{f}_{nh}(x)]}} \right) \xrightarrow{\mathfrak{D}} N(0,1),$$

*where $N(0,1)$ is the standard normal distribution.*

Asymptotic normality still holds if $nh^{2\beta+1} \to C$, where $C$ is some constant, but in this case the limit will not be distribution-free; it will depend on the unknown function $g$. We cannot select an optimal bandwidth to obtain (distribution-free) asymptotic normality, but this is also the case in ordinary kernel estimation. This fact comes from the trade-off between bias and variance, for the details, see the proof of the theorem. Now let us consider a different centering: $\hat{f}_{nh}(x) - \mathrm{E}[\hat{f}_{nh}(x)]$. Then the following theorem holds true.

**Theorem 2.2.** *Suppose that Conditions F, W, H and M hold, $g(x) \neq 0$ and $nh^{1+4\beta} \to 0$. Then we have*

$$\left( \frac{\hat{f}_{nh}(x) - \mathrm{E}[\hat{f}_{nh}(x)]}{\sqrt{\mathrm{Var}[\hat{f}_{nh}(x)]}} \right) \xrightarrow{\mathfrak{D}} N(0,1).$$

We see that, in this case, the additional condition on the bandwidth is weaker than the one in Theorem 2.1.

## 3. Simulation results and numerical aspects

In this section we present two simulations. They complement the asymptotic results of Theorems 2.1 and 2.2 and give some (although incomplete) indication of the finite sampling properties of the estimator.

In the first example, the true density $f$ is the standard normal density and $\lambda = 0.3$. The kernel we used is from Wand [18] and it has the rather complicated expression

$$w(t) = \frac{48t(t^2-15)\cos t - 144(2t^2-5)\sin t}{\pi t^7},$$

but its characteristic function looks much simpler and is given by

$$\phi_w(t) = (1-t^2)^3 1_{\{|t|<1\}}.$$

The estimator is based on 1000 observations and the bandwidth equals 0.14 (the bandwidth was selected by hand). To compute the estimator, we used the fast Fourier transform. The idea, which in spirit is close to the method for numerical evaluation of option prices proposed by Carr and Madan [3], is sketched as follows:

(i) Notice that whenever $\lambda < \log 2$, the distinguished logarithm in (1.2) reduces to the principal branch of the logarithm.

(ii) The main use of truncation in (1.3) is to prove asymptotic properties of the estimator and, in general, we do not need to use it in practice.

(iii) The computation of the empirical characteristic function can be significantly sped up by grouping the observations, the idea used to numerically evaluate ordinary kernel density estimators. However, we computed the empirical characteristic function directly, without grouping the observations. Notice that we do not use the values of the empirical characteristic function in its tails.

(iv) Notice that we can rewrite (1.2) as $f_{nh}(x) = f_{nh}^{(1)}(x) + f_{nh}^{(2)}(x)$, where

$$f_{nh}^{(1)}(x) = \frac{1}{2\pi\lambda} \int_0^\infty e^{-itx} \mathrm{Log}((e^\lambda - 1)\phi_{\mathrm{emp}}(t)\phi_w(ht) + 1)\,\mathrm{d}t,$$

$$f_{nh}^{(2)}(x) = \frac{1}{2\pi\lambda} \int_0^\infty e^{itx} \mathrm{Log}((e^\lambda - 1)\phi_{\mathrm{emp}}(-t)\phi_w(ht) + 1)\,\mathrm{d}t.$$

Using the trapezoid rule and setting $v_j = \eta(j-1)$, $f_{nh}^{(1)}(x)$ can be approximated by

$$f_{nh}^{(1)}(x) \approx \frac{1}{2\pi\lambda} \sum_{j=1}^N e^{-iv_j x} \psi(v_j) \eta.$$

Here we take $N$ to be some power of 2 and $\psi(v_j) = \mathrm{Log}((e^\lambda - 1)\phi_{g_{nh}}(v_j) + 1)$. The application of the Fast Fourier Transform to this sum will give us $N$ values of $f_{nh}^{(1)}$ and we employ a regular spacing size $\delta$, so that our values of $x$ are

$$x_u = -\frac{N\lambda}{2} + \delta(u-1),$$

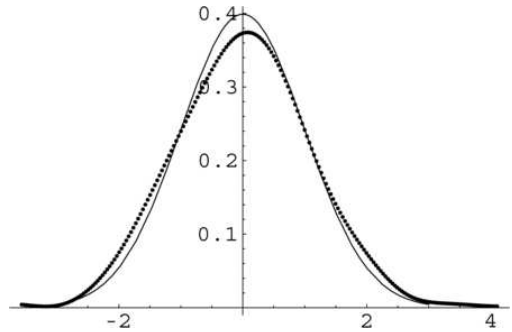

**Figure 1.** Estimation of a normal density.

where $u = 1, \ldots, N$. Thus we have

$$f_{nh}^{(1)}(x_u) \approx \frac{1}{2\pi\lambda} \sum_{j=1}^{N} e^{-i\delta\eta(j-1)(u-1)} e^{iv_j N\delta/2} \psi(v_j) \eta$$

for $u = 1, \ldots, N$. To apply the Fast Fourier Transform, we note that we must take $\delta\eta = 2\pi/N$. If we choose $\eta$ small to obtain a fine grid for integration, then we will obtain values of $f_{nh}^{(1)}$ at values of $x_u$ that are relatively seperate from each other. We would like, therefore, to obtain an accurate integration for larger values of $\eta$: to this end we incorporate Simpson weightings into our summation, that is,

$$f_{nh}^{(1)}(x_u) \approx \frac{1}{2\pi\lambda} \sum_{j=1}^{N} e^{-i(2\pi)/N(j-1)(u-1)} e^{iv_j N\delta/2} \psi(v_j) \frac{\eta}{3}(3 + (-1)^j - \delta_{j-1}),$$

where $\delta_j$ is a Kronecker function. Similar reasoning applies to $f_{nh}^{(2)}(x)$.

The result of this procedure for $N = 16\,384$ and $\eta = 0.01$ is given in Figure 1 (the estimate is represented by the bold dotted line).

In the second example we consider the case when $f$ is a mixture of two normal densities with means 0 and $3/2$ and variances 1 and $1/9$ with mixing probabilities $3/4$ and $1/4$, respectively. The estimator is based on 1000 observations and the bandwidth equals 0.1; the kernel is the same as in the first example. The result is given in Figure 2 (the estimate is plotted by the bold dotted line). Note that the estimator captures the bimodal character of the density $f$ in a quite satisfactory manner.

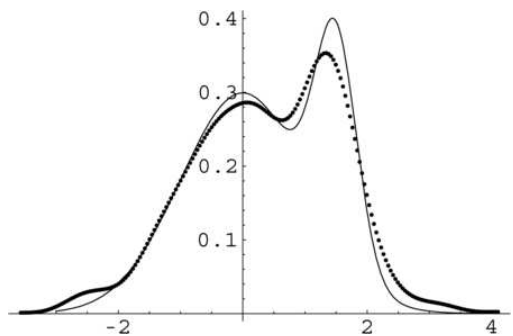

**Figure 2.** Estimation of a mixture of normal densities.

## 4. Proofs

**Proof of Lemma 2.1.** We have $|\phi_g(t)| \leq C|\phi_f(t)|$, which follows from the relationship

$$\phi_g(t) = \frac{1}{e^\lambda - 1}(e^{\lambda \phi_f(t)} - 1).$$

Indeed,

$$|e^{\lambda \phi_f(t)} - 1| = |-1 + 1 + \lambda \phi_f(t) + \cdots| \leq \lambda |\phi_f(t)| e^{\lambda |\phi_f(t)|} \leq \lambda e^\lambda |\phi_f(t)|$$

and $|\phi_g(t)| \leq C|\phi_f(t)|$ follows, where $C = \lambda e^\lambda (e^\lambda - 1)^{-1}$. This implies that $t^\beta \phi_g(t)$ is integrable. Furthermore,

$$g(x) = \sum_{n=1}^\infty f^{*n}(x) P(N = n | N > 0),$$

where $f^{*n}$ denotes the $n$-fold convolution of $f$. By Parseval's theorem,

$$\int_{-\infty}^\infty (g^{(l)}(x+t) - g^{(l)}(x))^2 \, dx = \int_{-\infty}^\infty |u^l \phi_g(u)|^2 |e^{itu} - 1|^2 \, du, \qquad (4.1)$$

where we used the fact that $|\phi_{g^{(l)}}(u)| = |u^l \phi_g(u)|$ (see Schwartz [13], pages 180–182). The latter is true because the derivatives of $g(x)$ up to order $l$ are integrable, which can be verified by direct computation employing formula (III, 2;8) of Schwartz [13]. From (4.1) it follows that

$$\int_{-\infty}^\infty (g^{(l)}(x+t) - g^{(l)}(x))^2 \, dx \leq \left(\frac{\lambda e^\lambda}{e^\lambda - 1}\right)^2 \int_{-\infty}^\infty |u^l \phi_f(u)|^2 |e^{itu} - 1|^2 \, du.$$

Applying Parseval's theorem to the right-hand side and recalling that $f$ belongs to $\mathcal{N}(\beta, L_2)$, we conclude that $g$ belongs to $\mathcal{N}(\beta, \lambda e^\lambda (e^\lambda - 1)^{-1} L_2)$. Now we will verify that $g \in \mathcal{H}(\beta, L_1)$. We have

$$g^{(l)}(x) = \sum_{n=1}^{\infty} f^{*(n-1)} * f^{(l)}(x) P(N = n | N > 0).$$

Using this expression, we get

$$|g^{(l)}(x+t) - g^{(l)}(x)|$$
$$= \left| \sum_{n=1}^{\infty} P(N = n | N > 0) \int_{-\infty}^{\infty} (f^{(l)}(x+t-u) - f^{(l)}(x-u)) f^{*(n-1)}(u) \, du \right|$$
$$\leq L_1 |t|^{\beta - l} \sum_{n=1}^{\infty} P(N = n | N > 0) \int_{-\infty}^{\infty} f^{*(n-1)}(u) \, du = L_1 |t|^{\beta - l}.$$

This completes the proof of the lemma. □

**Proof of Proposition 2.1.** We may write

$$b^w(n, h, x) = E[\hat{f}_{nh}(x) 1_{[J_n \leq \delta]} + \hat{f}_{nh}(x) 1_{[J_n > \delta]} - f(x) 1_{[J_n \leq \delta]}] - f(x) P(J_n > \delta),$$

where $\delta$ is any positive number and $J_n$ denotes the integrated squared error of the estimator $g_{nh}$. We have

$$|E[\hat{f}_{nh}(x) 1_{[J_n > \delta]}]| \leq M_n P(J_n > \delta).$$

This term is of order lower than $h^\beta$. To see this, recall the special form of $M_n$ and $h$, and apply the exponential bound to $P(J_n > \delta)$ that is valid for all $n$ sufficiently large (see Devroye [5], page 36, Remark 3). Also $f(x) P(J_n > \delta) = o(h^\beta)$.

Now we turn to

$$E[(\hat{f}_{nh}(x) - f(x)) 1_{[J_n \leq \delta]}].$$

By selecting $\delta$, we can achieve that $\phi_{g_{nh}}(t)$ is uniformly close to $\phi_g(t)$ on the set $\{J_n \leq \delta\}$. This is true because if $J_n \leq \delta$, then

$$|\phi_{g_{nh}}(t) - \phi_g(t)| = \left| \int_{-\infty}^{\infty} e^{-itx} (g_{nh}(x) - g(x)) \, dx \right| \leq J_n \leq \delta. \quad (4.2)$$

This in turn implies that for $\delta$ small (e.g., $\delta = e^{-\lambda}/2$), $(e^\lambda - 1)\phi_{g_{nh}}(t) + 1$ is bounded away from zero on the set $J_n \leq \delta$, because

$$|(e^\lambda - 1)\phi_g(t) + 1| = |e^{\lambda \phi_f(t)}| \geq e^{-\lambda}.$$

Therefore, the distinguished logarithm will be well defined on this set and $\log(|(e^\lambda - 1)\phi_{g_{nh}}(t) + 1)|$, that is, the real part of the distinguished logarithm $\mathrm{Log}((e^\lambda - 1)\phi_{g_{nh}}(t) + 1))$ will be bounded on $\{J_n \leq \delta\}$. The imaginary part of $\mathrm{Log}((e^\lambda - 1)\phi_{g_{nh}}(t) + 1)$ is also bounded. This holds true because, for $t$ sufficiently large, $(e^\lambda - 1)\phi_g(t) + 1$ is arbitrarily close to 1 and hence the argument of the distinguished logarithm $\mathrm{Log}((e^\lambda - 1)\phi_g(t) + 1)$ cannot circle around zero infinitely many times. To see the latter, we can argue as follows: there exists $t^*$ such that, for $t \geq t^*$, $(e^\lambda - 1)\phi_g(t) + 1$ does not make a turn around zero, because as $t \to \infty$, the function tends to $e^\lambda$. If we assume that $(e^\lambda - 1)\phi_{g_{nh}}(t) + 1$ in $[0, t^*]$ makes an infinite number of turns around zero, then its length on $[0, t^*]$ must also be infinite (because the curve stays away from zero at a positive distance). One can check that under given conditions on $w$, the latter is not true and, hence, also $(e^\lambda - 1)\phi_{g_{nh}}(t) + 1$ can make only a finite number of turns around zero.

Thus on the set $\{J_n \leq \delta\}$, the argument of $\mathrm{Log}((e^\lambda - 1)\phi_{g_{nh}}(t) + 1)$ will be bounded for $\delta$ small and, hence, on the set $\{J_n \leq \delta\}$ for large $n$ and small $\delta$, the truncation becomes unimportant and we have $\hat{f}_{nh}(x) = f_{nh}(x)$. Therefore,

$$\mathrm{E}[(f_{nh}(x) - f(x))1_{[J_n \leq \delta]}]$$
$$= \frac{1}{2\pi\lambda}\mathrm{E}\left[\left(\int_{-1/h}^{1/h} e^{-itx} \mathrm{Log}((e^\lambda - 1)\phi_{g_{nh}}(t) + 1)\,dt\right.\right.$$
$$\left.\left. - \int_{-1/h}^{1/h} e^{-itx} \mathrm{Log}((e^\lambda - 1)\phi_g(t) + 1)\,dt\right)1_{[J_n \leq \delta]}\right]$$
$$- \frac{1}{2\pi\lambda}\int_{-\infty}^{-1/h} e^{-itx}\phi_f(t)\,dt P(J_n \leq \delta) - \frac{1}{2\pi\lambda}\int_{1/h}^{\infty} e^{-itx}\phi_f(t)\,dt P(J_n \leq \delta).$$

The last two terms are of lower order than $h^\beta$. Indeed, we have, for example,

$$\lim_{h \to 0} \frac{1}{h^\beta}\left|\int_{1/h}^{\infty} e^{-itx}\phi_f(t)\,dt\right|$$
$$\leq \lim_{h \to 0} \frac{1}{h^\beta}\int_{1/h}^{\infty} |\phi_f(t)|\,dt \leq \lim_{h \to 0} \int_{1/h}^{\infty} t^\beta |\phi_f(t)|\,dt = o(1). \qquad (4.3)$$

Hence we need to study

$$\frac{1}{2\pi\lambda}\mathrm{E}\left[\left(\int_{-1/h}^{1/h} e^{-itx} \mathrm{Log}((e^\lambda - 1)\phi_{g_{nh}}(t) + 1)\,dt\right.\right.$$
$$\left.\left. - \int_{-1/h}^{1/h} e^{-itx} \mathrm{Log}((e^\lambda - 1)\phi_g(t) + 1)\,dt\right)1_{[J_n \leq \delta]}\right]$$
$$= \frac{1}{2\pi\lambda}\mathrm{E}\left[\int_{-1/h}^{1/h} e^{-itx} \mathrm{Log}(z_{nh}(t) + 1)\,dt 1_{[J_n \leq \delta]}\right], \qquad (4.4)$$

where

$$z_{nh}(t) = \frac{(e^\lambda - 1)(\phi_{g_{nh}}(t) - \phi_g(t))}{(e^\lambda - 1)\phi_g(t) + 1}.$$

Note that $z_{nh}$ is bounded. Rewrite (4.4) as

$$\frac{1}{2\pi\lambda}\mathrm{E}\left[\int_{-1/h}^{1/h} e^{-itx} z_{nh}(t)\,dt 1_{[J_n \leq \delta]}\right] + \frac{1}{2\pi\lambda}\mathrm{E}\left[\int_{-1/h}^{1/h} e^{-itx} R_{nh}(t)\,dt 1_{[J_n \leq \delta]}\right], \quad (4.5)$$

where

$$R_{nh}(t) = \mathrm{Log}(1 + z_{nh}(t)) - z_{nh}(t).$$

Consider the first term in (4.5). We claim that the omission of $1_{[J_n \leq \delta]}$ will result in an error of order lower than $h^\beta$. In fact,

$$\left|\mathrm{E}\left[\int_{-1/h}^{1/h} e^{-itx} z_{nh}(t)\,dt 1_{[J_n \leq \delta]}\right]\right|$$

$$\leq \left|\mathrm{E}\left[\int_{-1/h}^{1/h} e^{-itx} z_{nh}(t)\,dt\right]\right| - \left|\mathrm{E}\left[\int_{-1/h}^{1/h} e^{-itx} z_{nh}(t)\,dt 1_{[J_n > \delta]}\right]\right|.$$

The second term is bounded by $Ch^{-1}P(J_n > \delta)$, where $C$ is some constant, and this is of lower order than $h^\beta$ (recall the exponential bound of Devroye [5] on $P(J_n > \delta)$).

Using the fact that $\mathrm{E}[\phi_{\mathrm{emp}}(t)] = \phi_g(t)$, we obtain

$$\frac{1}{2\pi\lambda}\mathrm{E}\left[\int_{-1/h}^{1/h} e^{-itx} \frac{(e^\lambda - 1)(\phi_{\mathrm{emp}}(t)\phi_w(ht) - \phi_g(t))}{(e^\lambda - 1)\phi_g(t) + 1}\,dt\right]$$

$$= \frac{e^\lambda - 1}{2\pi\lambda}\int_{-1/h}^{1/h} e^{-itx} \frac{\phi_g(t)\phi_w(ht) - \phi_g(t)}{(e^\lambda - 1)\phi_g(t) + 1}\,dt$$

$$= \frac{e^\lambda - 1}{2\pi\lambda}\int_{-1/h}^{1/h} e^{-itx}(\phi_g(t)\phi_w(ht) - \phi_g(t))\,dt$$

$$+ \frac{e^\lambda - 1}{2\pi\lambda}\int_{-1/h}^{1/h} e^{-itx}(\phi_g(t)\phi_w(ht) - \phi_g(t))(e^{-\lambda\phi_f(t)} - 1)\,dt. \quad (4.6)$$

The first summand in the latter expression differs from the bias of the kernel estimator $g_{nh}(x)$ only by the absence of the term $-\int_{-\infty}^{-1/h}\phi_g(t)\,dt - \int_{1/h}^{\infty}\phi_g(t)\,dt$. This additional term is of lower order than $h^\beta$ (cf. (4.3)). Under Conditions W and F and due to Lemma 2.1, the bias of $g_{nh}(x)$ is of order $h^\beta$ (see Tsybakov [16], Proposition 1.2). As far as the second summand in (4.6) is concerned, it is dominated by

$$\lambda e^\lambda \frac{e^\lambda - 1}{2\pi\lambda}\int_{-1/h}^{1/h}|\phi_g(t)\phi_w(ht) - \phi_g(t)||\phi_f(t)|\,dt, \quad (4.7)$$

because
$$|e^{-\lambda \phi_f(t)} - 1| \leq \lambda e^\lambda |\phi_f(t)|.$$

Application of the Cauchy–Schwarz inequality to the integral in (4.7) yields that it is bounded from above by

$$\sqrt{\int_{-1/h}^{1/h} |\phi_g(t)\phi_w(ht) - \phi_g(t)|^2 \, dt} \sqrt{\int_{-1/h}^{1/h} |\phi_f(t)|^2 \, dt}.$$

The second factor in this expression is bounded uniformly in $h$ thanks to the fact that $\phi_f$ is integrable ($|\phi_f(t)|^2$ consequently is also integrable). As far as the first factor is concerned, by Parseval's theorem it is bounded by the integrated squared bias of the estimator $g_{nh}$,

$$\int_{-\infty}^{\infty} (g * w_h(x) - g(x))^2 \, dx,$$

where

$$w_h(x) = \frac{1}{h} w\left(\frac{x}{h}\right).$$

Because, under Conditions F and W, the integrated squared bias of $g_{nh}$ is of order $h^{2\beta}$ (see Tsybakov [16], Proposition 1.8), we conclude that (4.6) is of order $h^\beta$. This gives us the order of the leading term (4.6) in bias expansion.

Now we turn to the second term in (4.5). We have

$$\left| E\left[ \int_{-1/h}^{1/h} e^{-itx} R_{nh}(t) \, dt 1_{[J_n \leq \delta]} \right] \right| \leq E\left[ \int_{-1/h}^{1/h} |R_{nh}(t)| \, dt 1_{[J_n \leq \delta]} \right].$$

To deal with this term we will need the inequality

$$|\text{Log}(1 + z_{nh}(t)) - z_{nh}(t)| \leq |z_{nh}(t)|^2, \tag{4.8}$$

provided that $|z_{nh}(t)| < \frac{1}{2}$. This inequality follows from the inequality

$$|e^z - 1 - z| \leq z^2,$$

which is valid for $|z| < 1/2$ if we take $z = \text{Log}(1 + z_{nh}(t))$, because by choosing $n$ large enough and $\delta$ small, $J_n \leq \delta$ will entail $|z_{nh}(t)| < 1/2$; see (4.2). Using the inequality (4.8), we obtain

$$E\left[ \int_{-1/h}^{1/h} |R_{nh}(t)| \, dt 1_{[J_n \leq \delta]} \right] \leq E\left[ \int_{-\infty}^{\infty} |z_{nh}(t)|^2 \, dt \right]$$

$$\leq K E\left[ \int_{-\infty}^{\infty} |\phi_{\text{emp}}(t)\phi_w(ht) - \phi_g(t)|^2 \, dt \right]$$

$$= K \mathrm{E}\left[\int_{-\infty}^{\infty}(g_{nh}(t) - g(t))^2 \, \mathrm{d}t\right] = K \operatorname{MISE}_n(h), \quad (4.9)$$

where $K$ is a constant. Here we used the fact that $|(\mathrm{e}^\lambda - 1)\phi_g(t) + 1| = \mathrm{e}^{\lambda \phi_f(t)}$ is bounded from below and applied Parseval's identity. Using the bound on $\operatorname{MISE}_n(h)$ (see Tsybakov [16], page 21) and combining it with (4.6), we establish the desired result. □

**Proof of Proposition 2.2.** Throughout the proof we will frequently use the following version of the Cauchy–Schwarz inequality: if $\xi$ and $\eta$ are random variables, then $|\operatorname{Cov}[\xi, \eta]| \leq \sqrt{\operatorname{Var}[\xi]}\sqrt{\operatorname{Var}[\eta]}$ provided that the variances exist. Hence, if the variance of $\eta$ is negligible in comparison to that of $\xi$, then $\operatorname{Cov}[\xi, \eta]$ also will be negligible in comparison to $\operatorname{Var}[\xi]$ and, therefore, $\operatorname{Var}[\xi + \eta] \sim \operatorname{Var}[\xi]$; that is, the leading term of $\operatorname{Var}[\xi + \eta]$ is $\operatorname{Var}[\xi]$.

Now we turn to the proof of the proposition itself. We have

$$\operatorname{Var}[\hat{f}_{nh}(x)] = \operatorname{Var}[\hat{f}_{nh}(x)1_{[J_n \leq \delta]} + \hat{f}_{nh}(x)1_{[J_n > \delta]}].$$

The variance of $\hat{f}_{nh}(x)1_{[J_n > \delta]}$ is of lower order than $(nh)^{-1}$, because of the special form of $M_n = n^\alpha$, the exponential bound on $P(J_n > \delta)$ and the inequality

$$\operatorname{Var}[\hat{f}_{nh}(x)1_{[J_n > \delta]}] \leq \mathrm{E}[(\hat{f}_{nh}(x))^2 1_{[J_n > \delta]}] \leq M_n^2 P(J_n > \delta).$$

Therefore, it suffices to consider $\operatorname{Var}[\hat{f}_{nh}(x)1_{[J_n \leq \delta]}]$. We have

$$\operatorname{Var}[\hat{f}_{nh}(x)1_{[J_n \leq \delta]}] = \operatorname{Var}[\hat{f}_{nh}(x)1_{[J_n \leq \delta]} - f(x)]$$

and because again the variance of $f(x)1_{[J_n > \delta]}$ is of a lower order than $(nh)^{-1}$, we can consider $\operatorname{Var}[(\hat{f}_{nh}(x) - f(x))1_{[J_n \leq \delta]}]$ instead. As we have seen in the proof of Proposition 2.1, on the set $\{J_n \leq \delta\}$ for $n$ large and $\delta$ sufficiently small, $\hat{f}_{nh}(x) = f_{nh}(x)$ and the distinguished logarithm is well defined. Write

$$\operatorname{Var}[(f_{nh}(x) - f(x))1_{[J_n \leq \delta]}]$$

$$= \operatorname{Var}\left[\left(\frac{1}{2\pi\lambda}\int_{-1/h}^{1/h} \mathrm{e}^{-\mathrm{i}tx} z_{nh}(t)\, \mathrm{d}t + \frac{1}{2\pi\lambda}\int_{-1/h}^{1/h} \mathrm{e}^{-\mathrm{i}tx} R_{nh}(t)\, \mathrm{d}t\right.\right.$$

$$\left.\left. - \int_{1/h}^{\infty} \mathrm{e}^{-\mathrm{i}tx}\phi_f(t)\,\mathrm{d}t - \int_{-\infty}^{-1/h} \mathrm{e}^{-\mathrm{i}tx}\phi_f(t)\,\mathrm{d}t\right)1_{[J_n \leq \delta]}\right].$$

The variances of the last two terms are negligible. Indeed, we have, for example,

$$\operatorname{Var}\left[\left|\int_{1/h}^{\infty} \mathrm{e}^{-\mathrm{i}tx}\phi_f(t)\,\mathrm{d}t\right|1_{[J_n \leq \delta]}\right]$$

$$= \left|\int_{1/h}^{\infty} \mathrm{e}^{-\mathrm{i}tx}\phi_f(t)\,\mathrm{d}t\right|^2 \operatorname{Var}[1_{[J_n > \delta]}] \leq CP(J_n > \delta)$$

with some constant $C$.

Hence we have to deal with

$$\text{Var}\left[\left(\frac{1}{2\pi\lambda}\int_{-1/h}^{1/h} e^{-itx} z_{nh}(t)\,dt + \frac{1}{2\pi\lambda}\int_{-1/h}^{1/h} e^{-itx} R_{nh}(t)\,dt\right)1_{[J_n\leq\delta]}\right]$$
$$= \text{Var}[\text{I} + \text{II}]. \tag{4.10}$$

We show that II has a negligible variance compared to that of I. Indeed, using the bound (4.9) from the proof of Proposition 2.1,

$$nh\,\text{Var}\left[\left|\int_{-1/h}^{1/h} e^{-itx} R_{nh}(t)\,dt 1_{[J_n\leq\delta]}\right|\right]$$
$$\leq K^2 nh\,\text{E}[(\text{ISE}_n(h))^2]$$
$$= K^2 nh\,\text{Var}[\text{ISE}_n(h)] + K^2 nh(\text{MISE}_n(h))^2,$$

where $K$ is a constant. Due to the conditions $nh \to \infty$ and $nh^{1+4\beta} \to 0$, we see that $nh(\text{MISE}_n(h))^2$ tends to 0.

We deal with $nh\,\text{Var}[\text{ISE}_n(h)]$. Let us write the integrated squared error as

$$\text{ISE}_n(h) = \frac{1}{n^2 h}\sum_{j=1}^{n}\int_{-\infty}^{\infty}(w(t))^2\,dt + \frac{1}{n^2 h}\sum_{j\neq k} w*w\left(\frac{Z_j - Z_k}{h}\right)$$
$$- \frac{2}{nh}\sum_{j=1}^{n}\int_{-\infty}^{\infty} w\left(\frac{t - Z_j}{h}\right)g(t)\,dt + \int_{-\infty}^{\infty}(g(t))^2\,dt,$$

using that

$$\frac{1}{h}\int_{-\infty}^{\infty} w\left(\frac{t - Z_j}{h}\right)w\left(\frac{t - Z_k}{h}\right)dt = w*w\left(\frac{Z_j - Z_k}{h}\right)$$

because $w$ is symmetric. Here $w*w$ denotes the convolution of $w$ with itself. From this it follows that

$$nh\,\text{Var}[\text{ISE}_n(h)]$$
$$= \frac{1}{n^3 h}\text{Var}\left[\sum_{j\neq k} w*w\left(\frac{Z_j - Z_k}{h}\right) - 2n\sum_{j=1}^{n}\int_{-\infty}^{\infty} w\left(\frac{t - Z_j}{h}\right)g(t)\,dt\right]. \tag{4.11}$$

We study the variance of each term between the brackets in (4.11) separately. For the second term we have

$$\frac{1}{n^3 h}\text{Var}\left[2n\sum_{j=1}^{n}\int_{-\infty}^{\infty} w\left(\frac{t - Z_j}{h}\right)g(t)\,dt\right]$$

$$= \frac{4}{nh} \sum_{j=1}^{n} \text{Var}\left[\int_{-\infty}^{\infty} w\left(\frac{t - Z_j}{h}\right) g(t)\, \mathrm{d}t\right]$$

$$= \frac{4}{h} \text{Var}\left[\int_{-\infty}^{\infty} w\left(\frac{t - Z_1}{h}\right) g(t)\, \mathrm{d}t\right]$$

$$\leq \frac{4}{h} \mathrm{E}\left[\left(\int_{-\infty}^{\infty} w\left(\frac{t - Z_1}{h}\right) g(t)\, \mathrm{d}t\right)^2\right]. \tag{4.12}$$

Through a change of the integration variable it is easily seen that

$$\int_{-\infty}^{\infty} w\left(\frac{t - Z_1}{h}\right) g(t)\, \mathrm{d}t = h \int_{-\infty}^{\infty} w(u) g(uh + Z_1)\, \mathrm{d}u$$

$$\leq h A \int_{-\infty}^{\infty} |w(u)|\, \mathrm{d}u,$$

where we used the fact that $g$ is bounded. Hence (4.12) vanishes as $h \to 0$. Now we arrive at the computation of the variance of the first term between the brackets in (4.11). We have

$$\frac{1}{n^3 h} \text{Var}\left[\sum_{j \neq k} w * w\left(\frac{Z_j - Z_k}{h}\right)\right]$$

$$= \frac{4}{n^3 h} \text{Var}\left[\sum_{j < k} w * w\left(\frac{Z_j - Z_k}{h}\right)\right]$$

$$= \frac{4}{n^3 h} \sum_{i<j} \sum_{k<l} \text{Cov}\left[w * w\left(\frac{Z_i - Z_j}{h}\right), w * w\left(\frac{Z_k - Z_l}{h}\right)\right].$$

We have three possibilities:

1. $i, j, k, l$ are distinct. Then, because of the independence, the corresponding covariances are 0.
2. $i = k, j = l$. The number of such possibilities is of order $n^2$ and because the covariances in (4.12) are bounded (because the convolution $w * w$ is bounded), the sum of such terms will be of order $n^2$.
3. The last possibility is that three indices out of four are distinct, for example, $i = k, j \neq l$. The number of such terms is of order $n^3$. Thus we have to study the behaviour of for example,

$$\frac{1}{h} \text{Cov}\left[w * w\left(\frac{Z_i - Z_j}{h}\right), w * w\left(\frac{Z_i - Z_l}{h}\right)\right].$$

Writing out this covariance yields

$$\frac{1}{h}\operatorname{Cov}\left[w*w\left(\frac{Z_i-Z_j}{h}\right),w*w\left(\frac{Z_i-Z_l}{h}\right)\right]$$

$$=\frac{1}{h}\operatorname{E}\left[w*w\left(\frac{Z_i-Z_j}{h}\right)w*w\left(\frac{Z_i-Z_l}{h}\right)\right]$$

$$-\frac{1}{h}\left(\operatorname{E}\left[w*w\left(\frac{Z_i-Z_j}{h}\right)\right]\right)^2$$

$$\leq \frac{1}{h}\operatorname{E}\left[w*w\left(\frac{Z_i-Z_j}{h}\right)w*w\left(\frac{Z_i-Z_l}{h}\right)\right].$$

Note that because $w$ is bounded, therefore $w*w$ is also bounded and it is sufficient to study the behaviour of

$$\frac{1}{h}\operatorname{E}\left[\left|w*w\left(\frac{Z_i-Z_j}{h}\right)\right|\right]. \tag{4.13}$$

To do this, first note that $Z_i - Z_j$ has density

$$m(x) = \int_{-\infty}^{\infty} g(t-x)g(t)\,\mathrm{d}t.$$

Using the change of variable formula and Fubini's theorem, we see that (4.13) can be written as

$$\int_{-\infty}^{\infty}\frac{1}{h}\left|w*w\left(\frac{x}{h}\right)\right|m(x)\,\mathrm{d}x = \int_{-\infty}^{\infty}\int_{-\infty}^{\infty}\frac{1}{h}\left|w*w\left(\frac{x}{h}\right)\right|g(t-x)g(t)\,\mathrm{d}x\,\mathrm{d}t.$$

Due to the fact that $\lim_{|u|\to\infty}|w(u)|=0$ and applying the dominated convergence theorem, we conclude that this double integral converges to 0 as $h \to 0$. Hence (4.11) tends to zero. Thus $\operatorname{Var}[II]$ is indeed negligible in comparison to $\operatorname{Var}[I]$.

Now we need to study (cf. (4.10))

$$\operatorname{Var}\left[\frac{1}{2\pi\lambda}\int_{-1/h}^{1/h} \mathrm{e}^{-itx}z_{nh}(t)\,\mathrm{d}t\,1_{[J_n\leq\delta]}\right].$$

Once again, applying the by now standard argument, instead of $\int_{-1/h}^{1/h}$, we take $\int_{-\infty}^{\infty}$ and substitute $1_{[J_n\leq\delta]}$ with 1, because the error will be of a lower order than $(nh)^{-1}$. Furthermore,

$$\operatorname{Var}\left[\frac{1}{2\pi\lambda}\int_{-\infty}^{\infty}\mathrm{e}^{-itx}z_{nh}(t)\,\mathrm{d}t\right] = \operatorname{Var}[A_{nh}(x)+B_{nh}(x)],$$

where

$$A_{nh}(x) = \frac{e^\lambda - 1}{2\pi\lambda} \int_{-\infty}^{\infty} e^{-itx}(\phi_{\text{emp}}(t)\phi_w(ht) - \phi_g(t))\,dt = \frac{e^\lambda - 1}{\lambda}(g_{nh}(x) - g(x)),$$

$$B_{nh}(x) = \frac{e^\lambda - 1}{2\pi\lambda} \int_{-\infty}^{\infty} e^{-itx}(\phi_{\text{emp}}(t)\phi_w(ht) - \phi_g(t))(e^{-\lambda\phi_f(t)} - 1)\,dt.$$

For the variance of $g_{nh}(x)$ we have the expansion

$$\text{Var}[g_{nh}(x)] = \frac{1}{nh}g(x)\int_{-\infty}^{\infty}(w(t))^2\,dt + o\!\left(\frac{1}{nh}\right),$$

see Tsybakov ([16], Proposition 1.4).

We will show that the variance of $B_{nh}(x)$ is of a smaller order than $(nh)^{-1}$. Indeed,

$$nh\,\text{Var}[B_{nh}(x)] = nh\,\text{Var}\!\left[\frac{e^\lambda - 1}{2\pi\lambda}\int_{-\infty}^{\infty} e^{-itx}\!\left(\frac{1}{n}\sum_{j=1}^{n} e^{itZ_j}\phi_w(ht)\right)\!(e^{-\lambda\phi_f(t)} - 1)\,dt\right]$$

$$= \frac{(e^\lambda - 1)^2}{(2\pi\lambda)^2} h\,\text{Var}\!\left[\int_{-\infty}^{\infty} e^{-it(x-Z_1)}\phi_w(ht)(e^{-\lambda\phi_f(t)} - 1)\,dt\right].$$

Now note that

$$\left|\int_{-\infty}^{\infty} e^{-it(x-Z_1)}\phi_w(ht)(e^{-\lambda\phi_f(t)} - 1)\,dt\right| \leq \int_{-\infty}^{\infty}|e^{-\lambda\phi_f(t)} - 1|\,dt$$

and that the right-hand side is finite thanks to the fact that $\phi_f(t)$ is integrable. Because we have $\text{Var}[\xi] \leq K^2$, for a random variable $|\xi|$ bounded by a constant $K$, we conclude that $\text{Var}[B_{nh}(x)] = o(\frac{1}{nh})$.

By combining all the intermediate results, we see that the leading term of the $\text{Var}[\hat{f}_{nh}(x)]$ is

$$\frac{1}{nh}\frac{(e^\lambda - 1)^2}{\lambda^2} g(x)\int_{-\infty}^{\infty}(w(u))^2\,du$$

and that the other terms are of lower order than $(nh)^{-1}$. □

**Proof of Proposition 2.3.** The result follows immediately from the decomposition

$$\text{MSE}[\hat{f}_{nh}(x)] = \text{Var}[\hat{f}_{nh}(x)] + (b^w(n,h,x))^2$$

and Propositions 2.1 and 2.2. □

**Proof of Theorem 2.1.** The proof is based on repeated applications of Slutsky's theorem (see Serfling [14], Section 1.5.4); that is, we will show that we can separate a sequence that gives asymptotic normality from our normalized sum and show that the remainder term converges to zero in probability. Then Slutsky's theorem will imply that the normalized sum is itself asymptotically normal. Write

$$\frac{\hat{f}_{nh}(x) - f(x)}{\sqrt{\text{Var}[\hat{f}_{nh}(x)]}} = \frac{\hat{f}_{nh}(x) - f(x)}{\sqrt{\text{Var}[\hat{f}_{nh}(x)]}} 1_{[J_n \leq \delta]} + \frac{\hat{f}_{nh}(x) - f(x)}{\sqrt{\text{Var}[\hat{f}_{nh}(x)]}} 1_{[J_n > \delta]}. \qquad (4.14)$$

If we take $n$ large and $\delta$ small, then

$$\frac{\hat{f}_{nh}(x) - f(x)}{\sqrt{\text{Var}[\hat{f}_{nh}(x)]}} 1_{[J_n \leq \delta]} = \frac{f_{nh}(x) - f(x)}{\sqrt{\text{Var}[\hat{f}_{nh}(x)]}} 1_{[J_n \leq \delta]}.$$

We treat the first term in (4.14). We have

$$\frac{f_{nh}(x) - f(x)}{\sqrt{\text{Var}[\hat{f}_{nh}(x)]}} 1_{[J_n \leq \delta]}$$

$$= \frac{1}{\sqrt{\text{Var}[\hat{f}_{nh}(x)]}} 1_{[J_n \leq \delta]} \left( \frac{1}{2\pi\lambda} \int_{-1/h}^{1/h} e^{-itx} \text{Log}(1 + z_{nh}(t)) \, dt \right.$$

$$\left. - \frac{1}{2\pi} \int_{1/h}^{\infty} e^{-itx} \phi_f(t) \, dt - \frac{1}{2\pi} \int_{-\infty}^{-1/h} e^{-itx} \phi_f(t) \, dt \right)$$

$$= \frac{1}{\sqrt{\text{Var}[\hat{f}_{nh}(x)]}} 1_{[J_n \leq \delta]} \frac{1}{2\pi\lambda} \int_{-1/h}^{1/h} e^{-itx} \text{Log}(1 + z_{nh}(t)) \, dt$$

$$- \frac{1}{\sqrt{\text{Var}[\hat{f}_{nh}(x)]}} 1_{[J_n \leq \delta]} \frac{1}{2\pi} \int_{1/h}^{\infty} e^{-itx} \phi_f(t) \, dt$$

$$- \frac{1}{\sqrt{\text{Var}[\hat{f}_{nh}(x)]}} 1_{[J_n \leq \delta]} \frac{1}{2\pi} \int_{-\infty}^{-1/h} e^{-itx} \phi_f(t) \, dt. \qquad (4.15)$$

Let us denote the second and third expressions by I and II. We can write (4.15) as

$$\frac{1}{\sqrt{\text{Var}[\hat{f}_{nh}(x)]}} 1_{[J_n \leq \delta]} \frac{1}{2\pi\lambda} \int_{-1/h}^{1/h} e^{-itx} \text{Log}(1 + z_{nh}(t)) \, dt$$

$$- (\text{I} - \text{E}[\text{I}]) - (\text{II} - \text{E}[\text{II}]) - \text{E}[\text{I}] - \text{E}[\text{II}].$$

The second and third terms of this expression converge to zero in probability. This follows from the application of Chebyshev's inequality and the facts that

$$\text{Var}[1_{[J_n \leq \delta]}] = \text{Var}[1_{[J_n > \delta]}] \leq P(J_n > \delta) \sim e^{-Cn},$$

$$\mathrm{Var}[\hat{f}_{nh}(x)] \sim \frac{1}{nh}.$$

The application of Slutsky's theorem shows that we can neglect them. Now we take a further step and rewrite (4.15) as

$$\frac{1}{\sqrt{\mathrm{Var}[\hat{f}_{nh}(x)]}} 1_{[J_n \leq \delta]} \frac{1}{2\pi\lambda} \int_{-1/h}^{1/h} \mathrm{e}^{-\mathrm{i}tx} z_{nh}(t)\,\mathrm{d}t$$

$$+ \frac{1}{\sqrt{\mathrm{Var}[\hat{f}_{nh}(x)]}} 1_{[J_n \leq \delta]} \frac{1}{2\pi\lambda} \int_{-1/h}^{1/h} \mathrm{e}^{-\mathrm{i}tx} R_{nh}(t)\,\mathrm{d}t - \mathrm{E}[\mathrm{I} + \mathrm{II}].$$

Denote the second term in this expression by III. Rewrite the above expression as

$$\frac{1}{\sqrt{\mathrm{Var}[\hat{f}_{nh}(x)]}} 1_{[J_n \leq \delta]} \frac{1}{2\pi\lambda} \int_{-1/h}^{1/h} \mathrm{e}^{-\mathrm{i}tx} z_{nh}(t)\,\mathrm{d}t$$

$$+ (\mathrm{III} - \mathrm{E}[\mathrm{III}]) - \mathrm{E}[\mathrm{I} + \mathrm{II} - \mathrm{III}].$$

Again, $(\mathrm{III} - \mathrm{E}[\mathrm{III}])$ converges to zero in probability and, therefore, we can neglect it. After doing so, we rewrite the above expression as

$$\frac{1}{\sqrt{\mathrm{Var}[\hat{f}_{nh}(x)]}} 1_{[J_n \leq \delta]} \frac{1}{2\pi\lambda} \int_{-\infty}^{\infty} \mathrm{e}^{-\mathrm{i}tx} z_{nh}(t)\,\mathrm{d}t$$

$$- \frac{1}{\sqrt{\mathrm{Var}[\hat{f}_{nh}(x)]}} 1_{[J_n \leq \delta]} \frac{1}{2\pi\lambda} \int_{1/h}^{\infty} \mathrm{e}^{-\mathrm{i}tx} z_{nh}(t)\,\mathrm{d}t$$

$$- \frac{1}{\sqrt{\mathrm{Var}[\hat{f}_{nh}(x)]}} 1_{[J_n \leq \delta]} \frac{1}{2\pi\lambda} \int_{-\infty}^{-1/h} \mathrm{e}^{-\mathrm{i}tx} z_{nh}(t)\,\mathrm{d}t - \mathrm{E}[\mathrm{I} + \mathrm{II} - \mathrm{III}].$$

Denote the second and third terms in this expression by IV and V. Then we can write

$$\frac{1}{\sqrt{\mathrm{Var}[\hat{f}_{nh}(x)]}} 1_{[J_n \leq \delta]} \frac{1}{2\pi\lambda} \int_{-\infty}^{\infty} \mathrm{e}^{-\mathrm{i}tx} z_{nh}(t)\,\mathrm{d}t$$

$$- (\mathrm{IV} - \mathrm{E}[\mathrm{IV}]) - (\mathrm{V} - \mathrm{E}[\mathrm{V}]) - \mathrm{E}[\mathrm{I} + \mathrm{II} - \mathrm{III} + \mathrm{IV} + \mathrm{V}].$$

There is nothing random in IV and V except $1_{[J_n \leq \delta]}$. Due to Chebyshev's inequality, $(\mathrm{IV} - \mathrm{E}[\mathrm{IV}])$ and $(\mathrm{V} - \mathrm{E}[\mathrm{V}])$ converge to zero in probability and, therefore, can be neglected. We then have to deal with (recall the definition of $z_{nh}$)

$$\frac{1}{\sqrt{\mathrm{Var}[\hat{f}_{nh}(x)]}} 1_{[J_n \leq \delta]} \frac{\mathrm{e}^\lambda - \lambda}{\lambda} (g_{nh}(x) - g(x))$$

$$+ (\text{VI} - \text{E}[\text{VI}]) - \text{E}[\text{I} + \text{II} - \text{III} + \text{IV} + \text{V} - \text{VI}],$$

where

$$\text{VI} \equiv \frac{1}{\sqrt{\text{Var}[\hat{f}_{nh}(x)]}} 1_{[J_n \leq \delta]} \frac{1}{2\pi\lambda} \int_{-\infty}^{\infty} \mathrm{e}^{-\mathrm{i}tx} (\mathrm{e}^{\lambda} - 1)(\phi_{g_{nh}}(t) - \phi_g(t))(\mathrm{e}^{-\lambda\phi_f(t)} - 1) \, \mathrm{d}t.$$

The argument from the proof of Proposition 2.2 shows that the variance of VI converges to zero and, hence, by Chebyshev's inequality, $\text{VI} - \text{E}[\text{VI}]$ converges to zero in probability. Therefore, we can neglect it. Thus we have

$$\frac{1}{\sqrt{\text{Var}[\hat{f}_{nh}(x)]}} 1_{[J_n \leq \delta]} \frac{\mathrm{e}^{\lambda} - \lambda}{\lambda} (g_{nh}(x) - g(x))$$

$$- \text{E}[\text{I} + \text{II} - \text{III} + \text{IV} + \text{V} - \text{VI}].$$

Now rewrite this as

$$\frac{1}{\sqrt{\text{Var}[\hat{f}_{nh}(x)]}} 1_{[J_n \leq \delta]} \frac{\mathrm{e}^{\lambda} - \lambda}{\lambda} (g_{nh}(x) - \text{E}[g_{nh}(x)])$$

$$+ (\text{VII} - \text{E}[\text{VII}]) - \text{E}[\text{I} + \text{II} - \text{III} + \text{IV} + \text{V} - \text{VI} - \text{VII}], \quad (4.16)$$

where

$$\text{VII} \equiv \frac{1}{\sqrt{\text{Var}[\hat{f}_{nh}(x)]}} 1_{[J_n \leq \delta]} \frac{\mathrm{e}^{\lambda} - \lambda}{\lambda} (\text{E}[g_{nh}(x)] - g(x)).$$

Due to Chebyshev's inequality, $\text{VII} - \text{E}[\text{VII}]$ converges to zero in probability and, therefore, can be neglected. The asymptotically normal term stems from the first term in (4.16), because $1_{[J_n \leq \delta]} \to 1$ in probability and because

$$\left( \frac{g_{nh}(x) - \text{E}[g_{nh}(x)]}{\sqrt{\text{Var}[g_{nh}(x)]}} \right) \xrightarrow{\mathfrak{D}} N(0,1),$$

which can be verified along the lines of pages 61–62 of Prakasa Rao [12] by checking Lyapunov's condition. It is easy to see that

$$\text{E}[\text{I} + \text{II} - \text{III} + \text{IV} + \text{V} - \text{VI} - \text{VII}] = \text{E}\left[ \frac{f_{nh}(x) - f(x)}{\sqrt{\text{Var}[\hat{f}_{nh}(x)]}} 1_{[J_n \leq \delta]} \right].$$

Adding the second term in (4.14) to this expression results in

$$\frac{b^w(n,h,x)}{\sqrt{\text{Var}[\hat{f}_{nh}(x)]}} + \frac{\hat{f}_{nh}(x)1_{[J_n > \delta]} - \text{E}[\hat{f}_{nh}(x)1_{[J_n > \delta]}]}{\sqrt{\text{Var}[\hat{f}_{nh}(x)]}} - \frac{f(x)1_{[J_n > \delta]} - \text{E}[f(x)1_{[J_n > \delta]}]}{\sqrt{\text{Var}[\hat{f}_{nh}(x)]}}.$$

The first term goes to zero because we assume that $nh^{2\beta+1} \to 0$. Two other terms converge to zero in probability. Thus, thanks to Slutsky's theorem, these terms can be neglected and we establish the desired result. □

**Proof of Theorem 2.2.** Write

$$\hat{f}_{nh}(x) - \mathrm{E}[\hat{f}_{nh}(x)]$$
$$= (\hat{f}_{nh}(x) - f(x))1_{[J_n \leq \delta]} + (\hat{f}_{nh}(x) - f(x))1_{[J_n > \delta]} + (f(x) - \mathrm{E}[\hat{f}_{nh}(x)]).$$

Using the same type of arguments as in Theorem 2.1 (note that we will not need $nh^{2\beta+1} \to 0$, because the bias divided by the root of variance will be cancelled in intermediate computations), we see that we have to deal with

$$\frac{\mathrm{e}^\lambda - 1}{\lambda} \frac{g_{nh}(x) - \mathrm{E}[g_{nh}(x)]}{\sqrt{\mathrm{Var}[\hat{f}_{nh}(x)]}} - \frac{\hat{f}_{nh}(x) - f(x)}{\sqrt{\mathrm{Var}[\hat{f}_{nh}(x)]}} 1_{[J_n > \delta]}$$
$$- \frac{\mathrm{E}[(\hat{f}_{nh}(x) - f(x))1_{[J_n > \delta]}]}{\sqrt{\mathrm{Var}[\hat{f}_{nh}(x)]}}.$$

The first term gives asymptotic normality, while the last two terms tend to zero in probability. The application of Slutsky's theorem yields the desired result. □

# Acknowledgement

The research of the second author was financed by the Netherlands Organisation for the Advancement of Scientific Research (NWO).